\documentclass[11pt]{amsart}

\usepackage{amscd,amssymb,amsopn,amsmath,amsthm,mathrsfs,graphics,amsfonts,enumerate,verbatim,calc
}

\usepackage[all]{xy}

\usepackage[OT2,OT1]{fontenc}
\newcommand\cyr{%
\renewcommand\rmdefault{wncyr}%
\renewcommand\sfdefault{wncyss}%
\renewcommand\encodingdefault{OT2}%
\normalfont
\selectfont}
\DeclareTextFontCommand{\textcyr}{\cyr}

\usepackage{amssymb,amsmath}

\DeclareFontFamily{OT1}{rsfs}{}
\DeclareFontShape{OT1}{rsfs}{n}{it}{<-> rsfs10}{}
\DeclareMathAlphabet{\mathscr}{OT1}{rsfs}{n}{it}

\numberwithin{equation}{section}
\newtheorem{theorem}{Theorem}[section]
\newtheorem{lemma}[theorem]{Lemma}

\newtheorem{question}{Question}

\newtheorem*{maintheorem}{Main Theorem}

\theoremstyle{definition}
\newtheorem{definition}[theorem]{Definition}
\newtheorem{remark}[theorem]{Remark}
\theoremstyle{remark}

\newtheorem{acknowledgement}{Acknowledgement}

\renewcommand{\ker}{\operatorname{Ker}}

\newcommand{\fm}{\frak{m}}
\newcommand{\fp}{\frak{p}}
\newcommand{\fq}{\frak{q}}

\begin{document}
\title[Frobenius test exponents]{on the uniform bound of Frobenius test exponents}

\author[Pham Hung Quy]{Pham Hung Quy}
\address{Department of Mathematics, FPT University Hanoi, Vietnam}
\email{quyph@fe.edu.vn}

\thanks{
The author is partially supported by a fund of Vietnam National Foundation for Science
and Technology Development (NAFOSTED) under grant number
101.04-2017.10.}


\begin{abstract} In this paper we prove the existence of a uniform bound for Frobenius test exponents for parameter ideals of a local ring $(R, \fm)$ of prime characteristic in the following cases:
\begin{enumerate}
\item $R$ is generalized Cohen-Macaulay. Our proof is much simpler than the original proof of Huneke, Katzman, Sharp and Yao;
\item The Frobenius actions on all lower local cohomologies $H^i_{\fm}(R)$, $i < \dim R$, are nilpotent.
\end{enumerate}
\end{abstract}

\maketitle

\section{Introduction}
Let $R$ be a Noetherian commutative ring of prime characteristic $p>0$, and $I$ an ideal of $R$. The {\it Frobenius closure} of $I$ is $I^F = \{x \mid x^{p^e} \in I^{[p^e]} \text{ for some } e \ge 0\}$, where $I^{[p^e]} = (r^{p^e} \mid r \in I)$ is the $e$-th Frobenius power of $I$. It is hard to compute $I^F$. By the Noetherianness of $R$ there is an integer $e$, depending on $I$, such that $(I^F)^{[p^e]} = I^{[p^e]}$. We call the smallest number $e$ satisfying the condition the {\it Frobenius test exponent} of $I$, and denote it by $Fte(I)$. It is natural to expect the existence of a uniform number $e$, depending only on the ring $R$, such that for every ideal $I$ we have $(I^F)^{[p^e]} = I^{[p^e]}$, i.e. $Fte(I) \le e$ for every ideal $I$. If we have a positive answer to this question, then two conditions $x \in I^F$ and $x^{p^e} \in I^{[p^e]}$ are equivalent. This gives in particular a finite test for the Frobenius closure. However, Brenner \cite{B06} gave two-dimensional normal standard graded domains with no uniform bound for Frobenius test exponents of all ideals. In contrast, Katzman and Sharp \cite{KS06} showed the existence of a uniform bound for Frobenius test exponents if we restrict to the class of parameter ideals in a Cohen-Macaulay local ring. For any local ring $(R, \fm)$ the {\it Frobenius test exponent for parameter ideals} of $R$, denoted by $Fte(R)$, is the smallest integer $e$ such that $(\fq^F)^{[p^e]} = \fq^{[p^e]}$ for every parameter ideal $\fq$ of $R$, and $Fte(R) = \infty$ if there is no such integer. Katzman and Sharp asked whether $Fte(R) < \infty$ for any (equidimensional) local ring \cite[Introduction]{KS06}. Furthermore, the authors of \cite{HKSY06} confirmed the question for generalized Cohen-Macaulay local rings.\\

The main idea in \cite{HKSY06, KS06} is connecting the Frobenius test exponent for parameter ideals with an invariant defined by the Frobenius actions on local cohomology $H^i_{\fm}(R)$, $i \ge 0$, called {\it the Hartshorne-Speiser-Lyubeznik number} of $H^i_{\fm}(R)$. Recall that the Frobenius endomorphism $F: R \to R, x \mapsto x^p$ induces Frobenius actions on local cohomology $H^i_{\fm}(R)$ for all $i \ge 0$. Roughly speaking, {\it the Hartshorne-Speiser-Lyubeznik number} of $H^i_{\fm}(R)$, denoted by $HSL(H^i_{\fm}(R))$, is a nilpotency index of Frobenius action on $H^i_{\fm}(R)$ (see Section 2 for details). The Hartshorne-Speiser-Lyubeznik number of $R$ is $HSL(R) = \max \{HSL(H^i_{\fm}(R)) \mid i = 0, \ldots, \dim R\}$. It is proved that $Fte(R) = HSL(R)$ whenever $R$ is a Cohen-Macaulay local ring \cite{KS06}. In general, Huong and the author \cite{HQ18} showed that $Fte(R) \ge HSL(R)$ for any local ring.\\

In the case of generalized Cohen-Macaulay rings, besides using the  Hartshorne-Speiser-Lyubeznik number the authors of \cite{HKSY06} needed several concepts and techniques from commutative algebra, including unconditioned strong d-sequences, cohomological annihilators, and modules of generalized fractions. Their proof is highly technical and it is hard to apply for other cases. The motivation of the present paper is finding a simple proof for the main result of \cite{HKSY06}. Recently, the author studied the length $\ell_R(\fq^*/\fq)$, where $\fq^*$ denotes the tight closure of a parameter ideal $\fq$, in local rings $F$-rational on the functured spectrum \cite{Q18}. More precisely, if $(R, \fm)$ is a local ring of dimension $d$ that is $F$-rational on the punctured spectrum and is $F$-injective, then for every parameter ideal $\fq$ we have
$$\ell_R (\fq^* / \fq) = \sum_{i=0}^{d-1}\binom{d}{i} \ell_R (H^i_{\fm}(R)) + \ell_R (0^*_{H^d_{\fm}(R)}),$$
where $0^*_{H^d_{\fm}(R)}$ denotes the tight closure of the zero submodule of $H^d_{\fm}(R)$. The combinatorial coefficients in the formula indicate that we may study the tight closure and Frobenius closure of parameter ideals via quotient rings $R/(x)$, where $x$ is a parameter element. The Frobenius endomorphism $F:R/(x)\to R/(x)$ can be factored as $R/(x)\to R/(x^{p})\to R/(x)$ where the second map is the natural projection. We denote the first map by $F_R$. There are induced maps of local cohomology modules $F_R: H^i_\fm(R/(x))\to H^i_\fm(R/(x^{p}))$. We call $F_R$ the {\it relative Frobenius action on local cohomology} $H^i_\fm(R/(x))$ with respect to $R$ (see Section 2 for details). In \cite{PQ18} Polstra and the author used the relative Frobenius actions on local cohomology as the key ingredient to study {\it $F$-nilpotent} rings. Recall that a local ring $(R, \fm)$ of dimension $d$ is called $F$-nilpotent if the Frobenius actions on all $H^i_{\fm}(R)$, $i < \dim R$, and on $0^*_{H^d_{\fm}(R)}$ are nilpotent. The notion of relative Frobenius actions on local cohomology seems to be very useful for the question of Katzman and Sharp. Indeed, using this notion we not only give a simple proof for the main result of \cite{HKSY06} but also prove that $Fte(R)< \infty$ in a new case. The main result of this paper is as follows.

\begin{maintheorem}\label{main thm}
Let $(R, \fm)$ be a local ring of positive characteristic $p>0$. Then $Fte(R) < \infty$ in the following cases:
\begin{enumerate}
\item $R$ is generalized Cohen-Macaulay, i.e. $H^i_{\fm}(R)$ is finitely generated for all $i < \dim R$;
\item The Frobenius actions on all lower local cohomologies $H^i_{\fm}(R)$, $i < \dim R$, are nilpotent.
\end{enumerate}
\end{maintheorem}
In the next section we cover the basic notions and background material relevant to the results of last section. We will prove the Main Theorem in the last section.
\begin{acknowledgement} The author obtained the proof of the first part of the Main Theorem after discussions with Linquan Ma but he did not confirm the authorship of this paper. The second part of the Main Theorem was done during the joint work with Thomas Polstra on $F$-nilpotent rings. The notion of relative Frobenius actions on local cohomology became clear to the author via that joint work. The author is deeply grateful to the referee for careful reading of the paper and valuable suggestions.
\end{acknowledgement}

\section{Prelimitaries}
\subsection{Frobenius closure of ideals}
 Let $R$ be a Noetherian ring containing a field of prime characteristic $p>0$. Let $F:R \to R, x \mapsto x^p$ denote the Frobenius endomorphism. If we want to notationally distinguish the source and target of the $e$-th Frobenius endomorphism $F^e: R \xrightarrow{x \mapsto x^{p^e}} R$, we will use $F_*^e(R)$ to denote the target. $F_*^e(R)$ is an $R$-bimodule, which is the same as $R$ as
an abelian group and as a right $R$-module, that acquires its left $R$-module structure via the $e$-th Frobenius
endomorphism $F^e$. By definition the $e$-th Frobenius endomorphism $F^e: R \to F_*^e(R)$ sending $x$ to $F_*^e(x^{p^e}) = x \cdot F_*^e(1)$ is an $R$-homomorphism.
\begin{definition}[\cite{H96}] Let $I$ be an ideal of $R$ we define
\begin{enumerate}
\item The {\it $e$-th Frobenius power} of $I$ is $I^{[p^e]} = (x^{p^e} \mid x \in I)$.
  \item The {\it Frobenius closure} of $I$, $I^F = \{x \mid  x^{p^e} \in I^{[p^e]} \text{ for some } e \ge 0\}$.
\end{enumerate}
\end{definition}
\begin{remark} An element $x \in I^F$ if it is contained in the kernel of the composition
 $$R \to R/I \otimes_R R \xrightarrow{\mathrm{id} \otimes F^e} R/I \otimes_R F_*^e(R)$$
 for some $e \ge 0$. Moreover, since $R$ is Noetherian, $I^F$ is finitely generated. Therefore there exists an integer $e_0$ such that
 $$I^F = \mathrm{Ker}(R \to R/I \otimes_R R \xrightarrow{\mathrm{id} \otimes F^e} R/I \otimes_R F_*^e(R))$$
for all $e \ge e_0$.
\end{remark}
As the above discussion for every ideal $I$ there is an integer $e$ (depending on $I$) such that $(I^F)^{[p^e]} = I^{[p^e]}$.  We call the smallest integer satisfying this condition the {\it Frobenius test exponent} of $I$ and denote it by $Fte(I)$. A problem of Katzman and Sharp \cite[Introduction]{KS06} asks in its strongest form: does there exist a number $e$, depending only on the ring $R$, such that, for every ideal $I$ we have $(I^F)^{[p^e]} = I^{[p^e]}$. A positive answer to
this question, together with the actual knowledge of a bound for $e$, would give an algorithm to compute the Frobenius closure $I^F$. We call this integer $e$ a {\it Frobenius test exponent} for the ring $R$. Unfortunately, Brenner \cite{B06} gave two-dimensional normal standard graded domains with no finite Frobenius test exponent. In contrast, Katzman and Sharp showed the existence of a uniform bound for  Frobenius test exponents if we restrict to the class of parameter ideals of a Cohen-Macaulay local ring. It leads to the following natural question.
\begin{question}[Katzman-Sharp]\label{Katzman Sharp} Let $(R, \fm)$ be an (excellent equidimensional) local ring of prime characteristic $p$. Then does there exist an integer $e$ such that for every parameter ideal $\fq$ of $R$ we have $(\fq^F)^{[p^e]} = \fq^{[p^e]}$?
\end{question}
We define the {\it Frobenius test exponent for parameter ideals} of $R$, $Fte(R)$, the smallest integer $e$ satisfying the above question and $Fte(R) = \infty$ if we have no such integer $e$. Question \ref{Katzman Sharp} has an affirmative answer when $R$ is generalized Cohen-Macaulay by \cite{HKSY06}, and is open in general.

\subsection{Frobenius actions on local cohomology} The Frobenius test exponent for parameter ideals is closely related with an invariant defined by the Frobenius actions on local cohomology. For any ideal $I = (x_1, \ldots, x_t)$ local cohomology $H^i_I(R)$ may be computed as the cohomology of the \v{C}ech complex
$$
0 \to R \to \bigoplus_{i=1}^t R_{x_i} \to \cdots \to R_{x_1 \ldots x_t} \to 0.
$$
The Frobenius endomorphism $F:R \to R$ and its localizations induce a natural Frobenius action on local cohomology $F:H^i_I(R) \to H^i_{I^{[p]}}(R) \cong H^i_{I}(R)$ for all $i \ge 0$. There is a very useful way of describing the top local cohomology. It can be given as the direct limit of Koszul cohomologies
$$
H^t_I(R) \cong\varinjlim_n R/(x_1^n, \ldots, x_t^n),
$$
with the map in the system $\varphi_{n, m}: R/(x_1^n, \ldots, x_t^n) \to R/(x_1^m, \ldots, x_t^m)$ is the multiplication by $(x_1 \ldots x_t)^{m - n}$ for all $m \ge n$. Then for each $\overline{a} \in H^t_{I}(R)$, which is the canonical image of some $a+(x_1^n, \ldots, x_t^n)$, we find that $F(\overline{a})$ is the canonical image of $a^p +(x_1^{pn}, \ldots, x_t^{pn})$.

Notice that if $(R, \fm)$ is a local ring then $H^i_{\fm}(R)$ is always Artinian. In general, let $A$ be an Artinian $R$-module with a Frobenius action $F: A \to A$. Then we define the {\it Frobenius closure} $0^F_A$ of the zero submodule of $A$ is the submodule of $A$ consisting all elements $z$ such that $F^e(z) = 0$ for some $e \ge 0$. $0^F_A$ is the nilpotent part of $A$ by the Frobenius action. By \cite[Proposition 1.11]{HS77} and \cite[Proposition 4.4]{L97} there exists a non-negative integer $e$ such that $0^F_A = \ker (A \overset{F^e}{\longrightarrow} A)$ for all $i \ge 0$ (see also \cite{Sh07}). The smallest of such integers is called the {\it Hartshorne-Speiser-Lyubeznik number} of $A$ and denoted by $HSL(A)$. We define the {\it Hartshorne-Speiser-Lyubeznik number} of a local ring $(R, \frak m)$ as follows
$$HSL(R): = \min \{ e \mid   0^F_{H^i_{\fm}(R)} =   \ker (H^i_{\fm}(R) \overset{F^e}{\longrightarrow} H^i_{\fm}(R)) \text{ for all } i = 0, \ldots, \dim R\}.$$
If $R$ is Cohen-Macaulay, then Katzman and Sharp \cite{KS06} showed that $Fte(R)$ is equal to $HSL(R)$. Shimomoto and the author \cite[Main Theorem A]{QS17} proved that if every parameter ideal of $R$ is Frobenius closed, i.e. $Fte(R) = 0$, then $R$ is $F$-injective, i.e. $HSL(R) = 0$. Moreover we also constructed an $F$-injective ring with a non-Frobenius closed parameter ideal, i.e. $HSL(R) = 0$ but $Fte(R)>0$ (see \cite[Section 5]{QS17}). Recently, Huong and the author \cite[Theorem 3.4]{HQ18} showed that $Fte(R) \ge HSL(R)$ for any local ring of prime characteristic.

\subsection{Relative Frobenius actions on local cohomology} We now discuss the notion of relative Frobenius actions on local cohomology which plays the key role of this paper. This notion was introduced in \cite{PQ18}  in study $F$-nilpotent rings. We believe that this technique will be very useful for Question \ref{Katzman Sharp}. Let $K \subseteq I$ be ideals of $R$. The Frobenius endomorphism $F:R/K\to R/K$ can be factored as composition of two natural maps:
$$R/ K  \to R/K^{[p]} \twoheadrightarrow R/K,$$
where the second map is the natural projection map. We denote the first map by $F_R$: $F_R(a + K) = a^p + K^{[p]}$ for all $a \in R$. The homomorphism $F_R$ induces {\it the relative Frobenius actions on local cohomology} $F_R: H^i_I(R/K) \to H^i_I(R/K^{[p]})$ via \v{Cech} complexes. The {\it relative Frobenius closure of the zero submodule of $H^i_I(R/K)$ with respect to $R$} is defined as follows
$$0^{F_R}_{H^i_I(R/K)} = \{ \eta \mid F^e_R(\eta)=0 \in H^i_I(R/K^{[p^e]}) \text{ for some } e \gg 0 \}.$$
A local cohomology module $H^i_I(R/K)$ is called \emph{F-nilpotent with respect to $R$} if $H^i_I(R/K) = 0^{F_R}_{H^i_I(R/K)}$. Furthermore, if there exists an integer $e$ such that
$$0^{F_R}_{H^i_I(R/K)} = \mathrm{Ker}\big(H^i_I(R/K) \overset{F^e_R}{\longrightarrow} H^i_I(R/K^{[p^e]}) \big),$$
then we call the smallest of such integers the {\it Hartshorne-Speiser-Lyubeznik number of $H^i_I(R/K)$ with respect to $R$}, and denote it by $HSL_R(H^i_I(R/K))$. And convention that $HSL_R(H^i_I(R/K)) = \infty$ if we have no such integer.\\
%Let $(R, \fm)$ be a local ring of dimension $d$. We define the {\it Hartshorne-Speiser-Lyubeznik number of $R/K$ with respect to $R$}
%$$HSL_R(R/K):= \sup \{HSL_R(H^i_{\fm}(R/K)) \mid i = 0, \ldots, \dim R/K\}.$$
We will verify Question \ref{Katzman Sharp} in the case $H^i_{\fm}(R)$ is {\it $F$-nilpotent}, i.e. $0^{F}_{H^i_{\fm}(R)} =  H^i_{\fm}(R)$, for all $i < \dim R$. This class of rings was studied in both algebraic side and geometrical side recently \cite{PQ18, ST17}. We need the following result from \cite[Theorems 4.2 and 4.4]{PQ18}.
\begin{lemma}\label{nilpotent deform} Let $(R, \fm)$ be a local ring of dimension $d$ and of prime characteristic $p>0$ such that $H^i_{\fm}(R)$ is $F$-nilpotent for all $i < d$. Then for every filter regular sequence $x_1, \ldots, x_i$, $i< d$, we have $H^j_{\fm}(R/(x_1, \ldots, x_i))$ is $F$-nilpotent with respect to $R$ for all $j < d-i$.
\end{lemma}
\subsection{Generalized Cohen-Macaylay rings}  We recall the cohomological definition of generalized Cohen-Macaulay rings (see \cite{Tr86}, \cite[Section 4]{QS17}).
\begin{definition} Let $(R, \frak m)$ be a Noetherian local ring of dimension $d$, and $\fq = (x_1, \ldots, x_d)$ a parameter ideal of $R$. Then
\begin{enumerate}
\item The ring $R$ is called {\it generalized Cohen-Macaulay} if $H^i_{\fm}(R)$ is finitely generated for all $i<d$.
\item The parameter ideal $\fq$ is called {\it standard} if for all $i+j < d$ we have $\fq H^j_{\fm}(R/(x_1, \ldots, x_i)) = 0.$
\end{enumerate}
\end{definition}
The following are well-known facts of generalized Cohen-Macaulay rings.
\begin{remark}\label{R3.2} Let $R$ be a generalized Cohen-Macaulay local ring. Then
\begin{enumerate}
\item Every system of parameters $x_1, \ldots, x_d$ is a filter regular sequence, i.e. $x_i \notin \fp$ for all $\fp \in \mathrm{Ass}_R(R/(x_1, \ldots, x_{i-1})) \setminus \{\fm\}$ for all $i = 1, \ldots, d$.
\item There exists an integer $n_0$ such that for every parameter ideal $\fq = (x_1, \ldots, x_d)$ we have $(x_1^{n_0}, \ldots, x_d^{n_0})$ is standard.
\item Ma \cite{M15} showed that $HSL(R) =0$ if and only if $Fte(R) = 0$. Moreover, in these cases $R$ is Buchsbaum (see also \cite[Section 4]{QS17}).
\end{enumerate}
\end{remark}

\section{Proof of the main theorem}
We firstly take some observations. Let $\fq = (x_1, \ldots, x_d)$ be a parameter ideal of $R$. By the prime avoidance theorem we can assume henceforth that $x_1, \ldots, x_d$ is a filter regular sequence. For all $i = 0, \ldots, d$, we set $\fq_i = (x_1, \ldots, x_i)$. For all $i = 1, \ldots, d$, the multiplication map $R/\fq_{i-1} \xrightarrow{x_i} R/\fq_{i-1}$ induces the short exact sequence
$$0\rightarrow R/(\fq_{i-1}:x_i) \xrightarrow{x_i}R/\fq_{i-1} \rightarrow R/\fq_i\rightarrow 0.$$
Since $(\fq_{i-1} : x_i)/\fq_{i-1}$ has finite length, we have the induced exact sequence of local cohomology
$$\cdots \to H^j_{\fm}(R/\fq_{i-1}) \to H^j_{\fm}(R/\fq_i) \xrightarrow{\delta} H^{j+1}_{\fm}(R/\fq_{i-1}) \to \cdots $$
for all $j \le d-i$. For each $e \ge 0$, considering the commutative diagram
$$
\begin{CD}
R/\fq_{i-1} @>x_i>> R/\fq_{i-1}  \\
@VVF^{e}_RV @VVF^{e}_RV\\
R/\fq_{i-1}^{[p^e]}  @>x_i^{p^e}>> R/\fq_{i-1}^{[p^e]}. \\
\end{CD}
$$
It induces the following commutative diagram with exact rows
$$
\begin{CD}
0 @>>> R/(\fq_{i-1}:x_i) @>x_i>> R/\fq_{i-1} @>>> R/\fq_{i} @>>> 0   \\
@. @VV(F^{e}_R)'V @VVF^{e}_RV @VVF^{e}_RV\\
0 @>>> R/(\fq_{i-1}^{[p^e]}:x_i^{p^e})  @>x_i^{p^e}>> R/\fq_{i-1}^{[p^e]} @>>> R/\fq_{i}^{[p^e]} @>>> 0, \\
\end{CD}
$$
where the most left vertical map is the composition
$$R/(\fq_{i-1}:x_i) \overset{F^e_R}{ \longrightarrow } R/(\fq_{i-1}:x_i)^{[p^e]} \twoheadrightarrow R/(\fq_{i-1}^{[p^e]}:x_i^{p^e}). $$
Since $x_i^{p^e}, \ldots, x_i^{p^e}$ is a filter regular sequence, $(\fq_{i-1}^{[p^e]}:x_i^{p^e})/\fq_{i-1}^{[p^e]}$ has finite length. Thus
$$ H^{j+1}_{\fm}(R/(\fq_{i-1}^{[p^e]}:x_i^{p^e}) ) \cong  H^{j+1}_{\fm}(R/(\fq_{i-1}:x_i)^{[p^e]} ) \cong H^{j+1}_{\fm}(R/\fq_{i-1}^{[p^e]})$$
for all $j \ge 0$. Therefore for all $j \le d-i$ we have the following commutative diagram
$$
\begin{CD}
\cdots @>>> H^j_{\fm}(R/\fq_{i-1}) @>>> H^j_{\fm}(R/\fq_i)  @>\delta>>  H^{j+1}_{\fm}(R/\fq_{i-1}) @>>> \cdots \\
@. @VVF^{e}_RV @VVF^{e}_RV @VVF^{e}_RV (\star)\\
\cdots @>>> H^j_{\fm}(R/\fq_{i-1}^{[p^e]}) @>>> H^j_{\fm}(R/\fq_i^{[p^e]})  @>>>  H^{j+1}_{\fm}(R/\fq_{i-1}^{[p^e]})@>>> \cdots. \\
\end{CD}
$$
We are ready to prove the main result of this paper.
\begin{proof}[Proof of the Main Theorem] Let $\fq =(x_1, \ldots, x_d)$ be a parameter ideal of $R$, and $e$ an integer. Notice that $(\fq^F)^{[p^{e}]} \subseteq (\fq^{[p^{e}]})^F$, so we have $Fte(\fq) \le Fte(\fq^{[p^{e}]}) + e$. Moreover, we can see that $Fte(\fq) = HSL_R(H^0_{\fm}(R/\fq))$. Therefore it is enough to prove that in one of the following cases
\begin{enumerate}
\item $R$ is a generalized Cohen-Macaulay ring and $\fq = (x_1, \ldots, x_d)$ is a standard parameter ideal (see Remark \ref{R3.2} (2));
\item $H^i_{\fm}(R)$ is $F$-nilpotent for all $i<d$, and $\fq = (x_1, \ldots, x_d)$ is a parameter ideal
\end{enumerate}
we have
$$HSL_R(H^0_{\fm}(R/\fq)) \le \sum_{k=0}^d \binom{d}{k} HSL(H^k_{\fm}(R)).$$
In fact, we will prove a stronger fact that for all $i + j  \le d$ we have
$$HSL_R(H^j_{\fm}(R/(x_1, \ldots, x_i))) \le \sum_{k = j}^{i+j} \binom{i}{k-j} HSL(H^k_{\fm}(R)).$$
We will proceed by induction on $i$. The case $i = 0$ is trivial. Suppose $i>0$ and the assertion holds true for $i-1$. Let $\fq_i = (x_1, \ldots, x_i)$ for all $i = 0, \ldots, d$. Notice that the rows of diagram $(\star)$ are left exact sequences in the first case by the properties of standard parameter ideals.\\
For each $j \le d-i$ we set $ e_{j+1}: = HSL_R( H^{j+1}_{\fm}(R/\fq_{i-1}))$ and $e'_j = HSL_R(H^j_{\fm}(R/\fq_{i-1}^{[p^{e_{j+1}}]}) )$. Both $e_{j+1}$ and $e'_j$ are finite by the inductive hypothesis. We claim that $HSL_R(H^j_{\fm}(R/\fq_i)) \le e'_j + e_{j+1}$. Indeed, we consider the following commutative diagram
$$
\begin{CD}
H^j_{\fm}(R/\fq_{i-1}) @>>> H^j_{\fm}(R/\fq_i)  @>\delta>>  H^{j+1}_{\fm}(R/\fq_{i-1}) \\
@VVF^{e_{j+1}}_RV @VVF^{e_{j+1}}_RV @VVF^{e_{j+1}}_RV \\
 H^j_{\fm}(R/\fq_{i-1}^{[p^{e_{j+1}}]}) @>\alpha>> H^j_{\fm}(R/\fq_i^{[p^{e_{j+1}}]})  @>\beta>>  H^{j+1}_{\fm}(R/\fq_{i-1}^{[p^{e_{j+1}}]}) \\
 @VVF^{e'_{j}}_RV @VVF^{e'_{j}}_RV\\
 H^j_{\fm}(R/\fq_{i-1}^{[p^{e'_j + e_{j+1}}]}) @>>> H^j_{\fm}(R/\fq_i^{[p^{e'_j+ e_{j+1}}]}).\\
\end{CD}
$$
Taking any element $x \in 0^{F_R}_{ H^j_{\fm}(R/\fq_i) }$. We have $\delta(x) \in 0^{F_R}_{H^{j+1}_{\fm}(R/\fq_{i-1})}$, so $F^{e_{j+1}}_R(\delta(x)) = 0$. Hence we have $\beta(F^{e_{j+1}}_R(x) ) = 0$. Thus there is $y \in  H^j_{\fm}(R/\fq_{i-1}^{[p^{e_{j+1}}]})$ such that $\alpha(y) = F^{e_{j+1}}_R(x)$. Since  $x \in 0^{F_R}_{ H^j_{\fm}(R/\fq_i) }$, $F^{e_{j+1}}_R(x) \in 0^{F_R}_{ H^j_{\fm}(R/\fq_i^{[p^{e_{j+1}}]})}$. In the first case, the second map of each row is injective, so $y \in  0^{F_R}_{H^j_{\fm}(R/\fq_{i-1}^{[p^{e_{j+1}}]})}$.  In the second case, $H^j_{\fm}(R/\fq_{i-1}^{[p^{e_{j+1}}]})$ is $F$-nilpotent with respect to $R$ by Lemma \ref{nilpotent deform}. Therefore we always have $F^{e'_{j}}_R(y) = 0$. Now we can easily see that
$$F^{e'_j+ e_{j+1}}_R(x) = F^{e'_j}_R(F^{ e_{j+1}}_R(x)) = F^{e'_j}_R(\alpha(y)) =  0.$$
We have proved that  $HSL_R(H^j_{\fm}(R/\fq_i)) \le e'_j + e_{j+1}$. Applying the inductive hypothesis for both sequences $x_1, \ldots, x_{i-1}$ and $x_1^{p^{e_{j+1}}}, \ldots, x_{i-1}^{p^{e_{j+1}}}$ we have
$$e_{j+1} = HSL_R( H^{j+1}_{\fm}(R/\fq_{i-1})) \le \sum_{k = j+1}^{i+j} \binom{i-1}{k-j-1}HSL(H^k_{\fm}(R)),$$
and
$$e'_{j} = HSL_R(H^j_{\fm}(R/\fq_{i-1}^{[p^{e_{j+1}}]}) ) \le \sum_{k = j}^{i+j-1} \binom{i-1}{k-j} HSL(H^k_{\fm}(R)).$$
Therefore
\begin{eqnarray*}
HSL_R(H^j_{\fm}(R/\fq_i)) & \le &  \sum_{k = j}^{i+j-1} \binom{i-1}{k-j} HSL(H^k_{\fm}(R)) +  \sum_{k = j+1}^{i+j} \binom{i-1}{k-j-1}HSL(H^k_{\fm}(R))\\
& =& \sum_{k = j}^{i+j} \bigg(\binom{i-1}{k-j} + \binom{i-1}{k-j-1} \bigg) HSL(H^k_{\fm}(R))\\
& =& \sum_{k = j}^{i+j} \binom{i}{k-j} HSL(H^k_{\fm}(R)).
\end{eqnarray*}
The proof is complete.
\end{proof}

\end{document}